# CRITICAL EXPONENTS OF PLANAR GRADIENT PERCOLATION

By Pierre Nolin

*École Normale Supérieure and Université Paris-Sud*

We study gradient percolation for site percolation on the triangular lattice. This is a percolation model where the percolation probability depends linearly on the location of the site. We prove the results predicted by physicists for this model. More precisely, we describe the fluctuations of the interfaces around their (straight) scaling limits, and the expected and typical lengths of these interfaces. These results build on the recent results for critical percolation on this lattice by Smirnov, Lawler, Schramm and Werner, and on the scaling ideas developed by Kesten.

**1. Introduction.** We study some aspects of "gradient percolation." This is a model of inhomogeneous site percolation where the probability for each site to be occupied varies along some fixed direction, for instance, the $y$-axis in the plane. This model is described in the percolation survey paper [2]. It was introduced by physicists (see [10, 20]) to model phenomena like diffusion or chemical etching, for example, to describe the interface created by the welding of two pieces of metal. They argued that this is a case where some aspects of critical percolation can be empirically observed without fine-tuning any parameter (this is a self-critical system—the critical phenomenon appears spontaneously).

Since we will restrict our study to the triangular planar lattice, let us briefly describe the model in this particular case, even if it makes sense for other lattices and dimensions.

Recall that if one colors each cell of a honeycomb lattice independently in black or white with respective probability $p$ and $1-p$, then when $p > 1/2$, there is an infinite connected component of black cells and when $p < 1/2$, there is an infinite connected component of white cells. The value $1/2$ is called the *critical probability* of this homogeneous percolation model (and percolation at $p = 1/2$ is called *critical percolation*). This model has recently received a lot of attention, which has led to, among other things, the computation of "critical exponents" (see [23]).

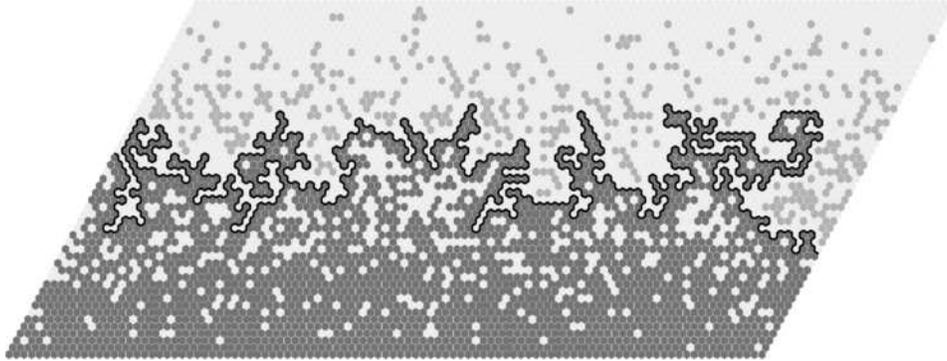

Fig. 1. *A simulation for $N = 50$ and $\ell_N = 100$.*

Suppose, now, that a large integer $N$ is given. We consider an inhomogeneous percolation model where each cell $z$ is colored in black and white independently, but with a probability that depends on $z$ (see Figure 1). More precisely, a cell $z$ with $y$-coordinate equal to $y(z) \in [0, N]$ will be colored in black with probability $p(z) = y/N$ [when $y < 0$, we take $p(z) = 0$ and when $y > N$, we take $p(z) = 1$]. It is then easy to see that there is almost surely a (unique) infinite black connected component (that contains the half-space $\{y > N\}$) and a unique white connected component (containing the half-space $\{y < 0\}$). Furthermore, the upper boundary of this white cluster and the lower boundary of the black cluster coincide. This separating curve is called the "*percolation front*."

Intuitively, it is quite clear what happens when $N$ is large. The percolation front will tend to be localized near the line $\{y/N = 1/2\}$. Furthermore, since at this level, the percolation is close to critical, the fine structure of the percolation front will be described in terms of critical percolation and its critical exponents.

Let us now briefly describe the main results of the present paper. Suppose that we restrict ourselves to a strip of height $N$ and of length $\ell_N$. Then, as $N \to \infty$ (provided $\ell_N$ also goes to infinity—not too fast but not too slowly— e.g., $\ell_N = N$ is acceptable for what follows), one can, with high probability, still define "the" percolation front separating the two "giant components". We call $\mathcal{R}_N$ the front (a precise definition will be given in Section 3) and $T_N$ its length (i.e., number of steps). We shall see that for each positive $\delta$, when $N \to \infty$, and with high probability:

- the front will remain in the strip of width $N^{4/7+\delta}$ near the line $\{y = N/2\}$;

- the front will not remain in the strip of width $N^{4/7-\delta}$ near the line $\{y = N/2\}$;
- the expected length $t_N = \mathbb{E}[T_N]$ of the front satisfies $N^{3/7-\delta}\ell_N \leq t_N \leq N^{3/7+\delta}\ell_N$;
- $T_N$ is close to its expected value, that is, $T_N/t_N$ is close to 1.

The localization of the front near the line $\{y = N/2\}$ and the values of the two exponents (4/7 for the width and 3/7 for the length) had been predicted by Sapoval et al. in [10, 20].

Our proofs build on the following mathematical results and ideas: Kesten's scaling relations [14], Smirnov's conformal invariance result [22] and the computation of the critical exponents for $SLE_6$ by Lawler, Schramm and Werner [15, 16] (Schramm Loewner Evolution, or SLE, was introduced by Schramm in [21]). More precisely, we shall further develop some ideas introduced in the first paper [14] and directly use results of [23], where Smirnov and Werner showed how to derive the critical exponents for percolation from those for $SLE_6$. We will not mention $SLE_6$ in our proofs, but we would like to stress that it plays a crucial role in the derivation of the results of [23] that we will use.

In fact, gradient percolation is rather easy to simulate. It has turned out to be an efficient practical tool to obtain numerical estimates for the critical probability of percolation on various lattices (by using the mean height of the front in order to approximate $p_c$—see, e.g., [19, 25]), for instance, the square lattice, and it has also been one of the first ways to get numerical evidence for values of the critical exponents of standard percolation (which then supported the conjectures based on Coulomb gas and conformal field theory).

Let us stress that the anisotropy of this model yields that despite the fact that the front converges to a straight line in the fine mesh limit (i.e., a curve of dimension one), its length for a lattice approximation of mesh-size $\delta$ (in a rhombus) behaves roughly like $\delta^{-3/7}$ (i.e., it has $\delta^{-10/7}$ steps). One can expect to observe a nontrivial limit—of fractal dimension 7/4—with an appropriate scaling (in $N^{4/7}$) of the axes, but the critical exponents obtained do not correspond directly to a fractal dimension for the limiting object.

**2. Homogeneous percolation preliminaries.** We recall in this section some known facts concerning percolation (not gradient percolation) that we will use later.

2.1. *Setting.* The setting in this paper will be site percolation in two dimensions on the triangular lattice. We will represent it, as usual, as a random (black or white) coloring of the faces of the hexagonal lattice. It has been proven by Kesten [12] and is now a classical fact that for this model, the critical probability is $p_c = 1/2$ and there is almost surely no infinite cluster when the percolation parameter is taken to be $1/2$. The reason why we focus here on this lattice is that it is (at present) the only one for which conformal invariance in the scaling limit has been proven (Smirnov [22]). Conformal invariance will not be used directly in the following, but this property, combined with the study of SLE by Lawler, Schramm and Werner [15, 16], allows the rigorous computation of the so-called "critical exponents" that will be instrumental in our considerations. On other regular lattices, like the square lattice $\mathbb{Z}^2$, these exponents are expected to be the same, but no rigorous derivation yet exists. However, some inequalities can be proven which would imply weaker, but nevertheless interesting, statements. We will not concern ourselves with this in the present paper.

The percolation parameter will be denoted by $p$: each site is *occupied* (or black) with probability $p$ and *vacant* (white) otherwise, independently of each other. The corresponding probability measure on the set of configurations will be referred to as $\mathbb{P}_p$, and $\mathbb{E}_p$ will denote the expectation.

We will use oblique coordinates, with the origin at 0 and the basis given by 1 and $e^{i\pi/3}$ (in complex notation), as shown on Figure 2. The parallelogram $R$ with corners $a_j + b_k e^{i\pi/3}$ ($j, k = 1, 2$) will thus be denoted by $[a_1, a_2] \times [b_1, b_2]$. Its interior and boundary will be denoted by $\overset{\circ}{R}$ and $\partial R$, respectively. We will consider $\|z\|$, the infinity norm of a site $z$, as measured with respect to the chosen coordinates, and $d(z, z') := \|z - z'\|$, the associated distance. For a site $z = (z_1, z_2)$, we will often use the rhombus $S_n(z) := [z_1 - n, z_1 + n] \times [z_2 - n, z_2 + n]$, which is the set of sites at a distance at most $n$ from $z$ [its interior, $\overset{\circ}{S}_n(z)$, and its boundary, $\partial S_n(z)$, consisting of the sites at a distance strictly less than $n$, resp. exactly $n$]. We will refer to $S_n(0)$ simply as $S_n$ and call it the "box of size $n$."

For two positive functions $f$ and $g$, $f \asymp g$ means that there exist two positive and finite constants $C_1$ and $C_2$ such that $C_1 g \leq f \leq C_2 g$ (so that their ratio is bounded away from 0 and $+\infty$) and $f \approx g$ means that $\log f / \log g \to 1$ (when $p \to 1/2$ or when $n \to \infty$, which will be clear from the context).

We now recall some relevant results on critical percolation that we will use in the present paper.

2.2. *Arm exponents.* We first briefly recall some facts concerning critical exponents for the existence of a certain number of "arms." These exponents describe the asymptotic behavior of the probability of certain exceptional events.

Let us consider a fixed integer $j \geq 2$. For any positive integers, we define the event $A^j(m, n)$ that there exist $j$ disjoint monochromatic paths from $\partial S_m$ to $\partial S_n$ that are not all of the same color (each path is either completely black/occupied or completely white/vacant, and there is at least one vacant path and one occupied path). As observed in [1, 23], we could also prescribe the cyclic order of the paths without changing the results that we state below.

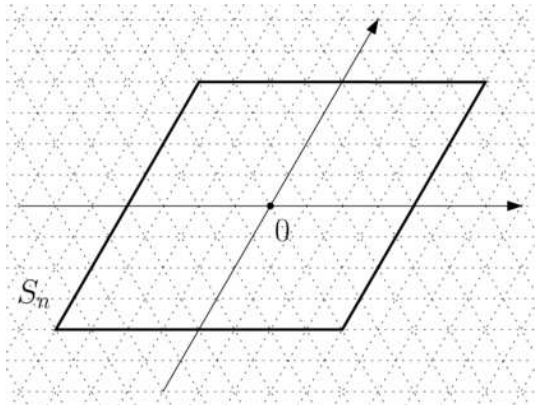

Fig. 2. *The triangular lattice, its associated basis and $S_n$.*

Combining the property of conformal invariance in the scaling limit (see [5, 22]) with the study of SLE made by Lawler, Schramm and Werner [15, 16], the following has been proven.

PROPOSITION 1 ([23]). *For all fixed $j \geq 2$, $m \geq j$,*

$$\mathbb{P}_{1/2}(A^j(m,n)) \approx n^{-(j^2-1)/12} \tag{1}$$

*when $n \to \infty$.*

In fact, we will use this result only for $j = 2, 3$ and 4. Let us simply remark that to derive these three exponents, it might be possible to bypass the use of the rather involved results of Camia and Newman [5], as these are exponents related to outer boundaries of clusters so that "universality"-based ideas might just be enough.

The value of the related "one-arm" exponent is 5/48 and has been derived in [17], but we shall not need it in the present paper.

2.3. *Behavior near criticality.* In the seminal paper [14], Kesten showed how the behavior of percolation at its critical point was related to the asymptotic behavior of percolation near its critical point and derived the so-called scaling relations that link some of the previous arm exponents to other critical exponents describing the behavior of connectivity probabilities near $p = p_c$.

A key idea in that article is to use a certain "characteristic length" $L(p)$ defined in terms of crossing probabilities (sometimes referred to as "sponge-crossing probabilities"). This idea (or adaptations of it) was also used in papers concerning finite-size scaling, for example, [3, 4, 8, 9].

Let us introduce some more notation. We will denote by $\mathcal{C}_H([a_1,a_2] \times [b_1,b_2])$ (resp. $\mathcal{C}_V([a_1,a_2] \times [b_1,b_2])$) the event that there exists a horizontal (resp. vertical) occupied crossing of the parallelogram $[a_1,a_2] \times [b_1,b_2]$, and by $\mathcal{C}_H^*, \mathcal{C}_V^*$, the same events with vacant crossings. However, for convenience, we will use a definition slightly different from the usual one. First, we decide to relax the condition on the boundary sites: a horizontal occupied crossing of a parallelogram $R$ will be a path connecting its left and right sides, all the sites of which *except its two extremities* are in $\overset{\circ}{R}$ and occupied. The events $\mathcal{C}_H(R)$ and $\mathcal{C}_H(R')$ will thus be independent if $R$ and $R'$ share one side. We also allow such a crossing to incorporate "dangling ends": it is not necessarily a self-avoiding path.

We have, for example,

(2) $\quad \mathbb{P}_p(\mathcal{C}_H([a_1,a_2] \times [b_1,b_2])) = 1 - \mathbb{P}_p(\mathcal{C}_V^*([a_1,a_2] \times [b_1,b_2]))$.

Now, consider rhombi $[0,n] \times [0,n]$. At $p=1/2$, $\mathbb{P}_p(\mathcal{C}_H([0,n] \times [0,n])) = 1/2$. When $p < 1/2$ (subcritical regime), this probability tends to 0 when $n$ goes to infinity, and it tends to 1 when $p > 1/2$ (supercritical regime).

We define a quantity that will roughly measure the scale up to which these crossing probabilities remain bounded away from 0 and 1. For each fixed $\varepsilon_0 > 0$, we define

(3) $\quad L(p, \varepsilon_0) = \begin{cases} \min\{n \text{ s.t. } \mathbb{P}_p(\mathcal{C}_H([0,n] \times [0,n])) \leq \varepsilon_0\}, & \text{when } p < 1/2, \\ \min\{n \text{ s.t. } \mathbb{P}_p(\mathcal{C}_H^*([0,n] \times [0,n])) \leq \varepsilon_0\}, & \text{when } p > 1/2. \end{cases}$

If we use the Russo–Seymour–Welsh (RSW) theory (see, e.g., [11, 13]), we see that for each $k \geq 1$, there exists some $\delta_k > 0$ (depending only on $\varepsilon_0$) such that

(4) $\quad \forall n \leq L(p) \quad \mathbb{P}_p(\mathcal{C}_H([0, kn] \times [0, n])) \geq \delta_k$.

For symmetry reasons, this bound is also valid for horizontal vacant crossings.

These estimates for crossing probabilities are then the basic building blocks on which many further considerations are built, for instance, the following one-arm probability estimate, which is one of the main results of Kesten's paper [14].

(5) $\quad \mathbb{P}_p[0 \rightsquigarrow \partial S_n] \asymp \mathbb{P}_{1/2}[0 \rightsquigarrow \partial S_n]$

for all $n \leq L(p)$ (uniformly in $p$). In Section 4.1 of the present paper, we will derive the analogous result for "two-arm probabilities."

This result is basically saying that when $n$ is not larger than $L(p)$, the situation can be compared to critical percolation. On the other hand, the definition of $L(p)$ shows that when $n > L(p)$, the picture starts to look like super/subcritical percolation. For instance, we have the following.

LEMMA 2 [Exponential decay with respect to $n/L(p)$]. *If $\varepsilon_0$ has been chosen sufficiently small, then there exists a constant $C > 0$ such that for all $n$ and all $p < 1/2$,*

(6) $$\mathbb{P}_p(\mathcal{C}_H([0,n] \times [0,n])) \leq C e^{-n/L(p)}.$$

Variants of this result are implicitly used or mentioned in Kesten's paper [14] and other papers on finite-size scaling. We now give its proof, as it just takes a couple of lines and clarifies matters.

PROOF OF LEMMA 2. Observe, first, that for all integer $n$,

(7) $$\mathbb{P}_p(\mathcal{C}_H([0,2n] \times [0,4n])) \leq C'[\mathbb{P}_p(\mathcal{C}_H([0,n] \times [0,2n]))]^2$$

with $C' = 10^2$ some universal constant. It suffices for that (see Figure 3) to divide the parallelogram $[0,2n] \times [0,4n]$ into four horizontal subparallelograms $[0,2n] \times [in, (i+1)n]$ ($i = 0, \ldots, 3$) and six vertical ones $[in, (i+1)n] \times [jn, (j+2)n]$ ($i = 0, 1$, $j = 0, \ldots, 2$). Indeed, consider a horizontal crossing of the big parallelogram: we can extract from it two pieces, one between its extremity on the left side and its first intersection with the vertical median $x = n$ and, in the same way, another one starting from the right side. These two subpaths both cross one of the subparallelograms "in the easy way." As they are disjoint by construction, the claim follows by using the van den Berg–Kesten inequality ([11, 24]).

We then obtain, by iterating,

(8) $$C'\mathbb{P}_p(\mathcal{C}_H([0, 2^k L(p)] \times [0, 2^{k+1} L(p)])) \leq (C'\varepsilon_1)^{2^k},$$

as soon as $\varepsilon_1 \geq \mathbb{P}_p(\mathcal{C}_H([0, L(p)] \times [0, 2L(p)]))$.

Recall that by definition, $\mathbb{P}_p(\mathcal{C}_H([0, L(p)] \times [0, L(p)])) \leq \varepsilon_0$. Consequently, the RSW theorem implies that for all fixed $\varepsilon_1 > 0$, if we take $\varepsilon_0$ sufficiently small, we automatically get (and independently of $p$) that

(9) $$\mathbb{P}_p(\mathcal{C}_H([0, L(p)] \times [0, 2L(p)])) \leq \varepsilon_1.$$

We now choose $\varepsilon_1 = 1/(e^2 C')$. For each integer $n \geq L(p)$, we can define $k = k(n)$ such that $2^k \leq n/L(p) < 2^{k+1}$ and then

$$\mathbb{P}_p(\mathcal{C}_H([0,n] \times [0,n])) \leq \mathbb{P}_p(\mathcal{C}_H([0, 2^k L(p)] \times [0, 2^{k+1} L(p)]))$$
$$\leq e^{-2^{k+1}}/C'$$
$$\leq e \times e^{-n/L(p)},$$

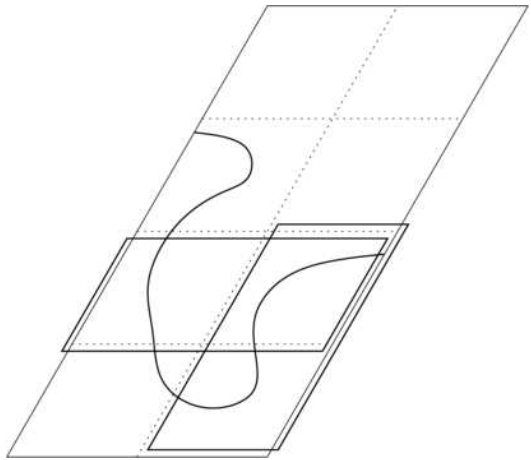

Fig. 3. *A horizontal crossing of the parallelogram* $[0, 2n] \times [0, 4n]$.

which is also valid for $n < L(p)$, due to the extra factor $e$. □

By counting the sites which are pivotal for the existence of a crossing, as outlined in Kesten's paper [14] (see the remark following Lemma 8 page 144–145, and also [18], for more details), we can prove that

$$(10) \qquad |p - 1/2|(L(p))^2 \mathbb{P}_{1/2}(A^4(1, L(p))) \asymp 1.$$

Combining this with the estimate of Proposition 1 for four arms leads to the following.

PROPOSITION 3 ([14, 23]). *When* $p \to 1/2$,

$$(11) \qquad L(p) \approx |p - 1/2|^{-4/3}.$$

REMARK 4. It has been shown in [14] that for any fixed $\varepsilon_1$ and $\varepsilon_2$ with $0 < \varepsilon_1, \varepsilon_2 \leq \varepsilon_0$,

$$(12) \qquad L(p, \varepsilon_1) \asymp L(p, \varepsilon_2).$$

Thus, the particular choice of $\varepsilon_0$ is not really important here (as long as it is sufficiently small).

## 3. Localization of the front.

3.1. *Framework of gradient percolation.* We will now define the model itself and fix some notation. The starting point is still site percolation on the triangular lattice. We first consider a strip $\mathcal{S}_N$ of finite width $2N$ (we assume it to be even for simplicity), centered around the $x$-axis, such that the $y$-coordinate varies between $-N$ and $N$ (see Figure 4). This strip may, for the moment, be unbounded in one or both directions and we will denote its length by $\ell_N$.

In this strip, we assume the parameter to decrease linearly according to $y$, that is, we choose it to be

(13) $$p(y) = 1/2 - y/2N.$$

With this choice, all the sites on the bottom edge $\mathcal{B}_N$ will be occupied ($p=1$) and all the sites on the top edge $\mathcal{T}_N$ vacant ($p=0$). The corresponding probability measure will be denoted by $\mathbb{P}$.

When we perform such a percolation, two opposite regions appear. At the bottom of $\mathcal{S}_N$, the parameter is close to 1, we are in a supercritical region and most occupied sites are connected to the bottom edge. On the contrary, we observe a big cluster of vacant sites on the top. The characteristic phenomenon of this model is the existence of a unique "front," a continuous interface touching both the occupied sites connected (by an occupied path) to the bottom of the strip and the vacant sites connected to the top (by a vacant path). Temporarily, we adopt the following definition.

DEFINITION 5. A *front* will be any path $\gamma_N$ on the dual hexagonal lattice that is bordered by an occupied (horizontal) crossing of $\mathcal{S}_N$ on one side and by a vacant (horizontal) crossing on the other side.

We will prove at the end of this section that under suitable hypotheses, with high probability, there is a unique such interface of this type. Note that when the length $\ell_N$ of the strip is finite, there is not necessarily a *unique* front. For instance, there is a positive probability to observe two horizontal crossings, one vacant and one occupied above it. If the strip is infinite, by independence of the different columns, there exists almost surely a column on which all the sites, except the highest one, are occupied. In that case, the front is unique and can be determined "dynamically" by starting two exploration processes from the top of this random column in the two possible directions (these processes follow the interface between black and white sites, making, at each step, a 60 degree turn to the left or right, according to the color of the hexagon they meet).

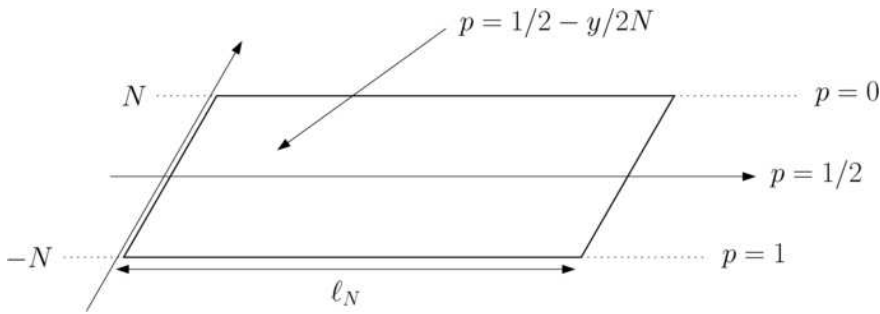

FIG. 4. *The strip where the percolation parameter $p$ varies.*

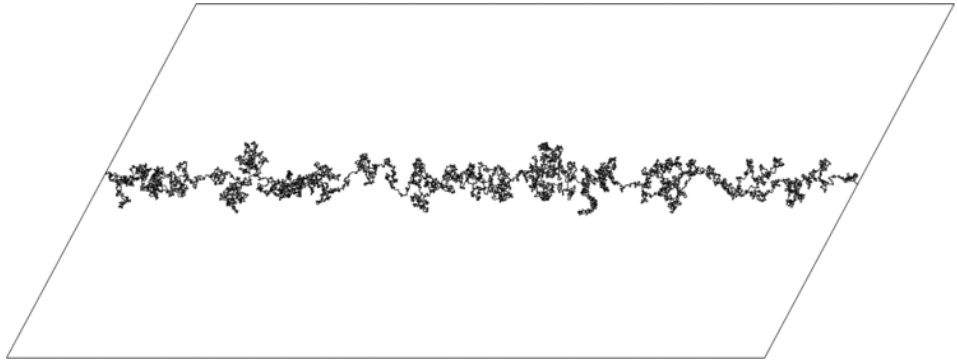

Fig. 5. *Localization of the front ($2N = 500$, $\ell_N = 1000$).*

However, for practical purposes, it will be more natural to consider strips of finite length. The argument above shows that on an infinite strip, exceptional events occur (the front touches the top and the bottom of the strip), and we may expect them not to happen when the length $\ell_N$ is not too large compared to $N$.

We now use the strip $\mathcal{S}_N = [0, \ell_N] \times [-N, N]$. We will often consider substrips of this big strip. For the sake of simplicity, the strip $[0, \ell_N] \times [-\lfloor N^\alpha \rfloor, \lfloor N^\alpha \rfloor]$ will be denoted by $[\pm N^\alpha]$.

3.2. *Localization.* We now study the convergence of the front when it is unique. Temporarily, we will consider, instead of "the" front, the highest horizontal crossing, and we will denote it by $\mathcal{R}_N$. Note that the sites just above $\mathcal{R}_N$ form a vacant horizontal crossing of $\mathcal{S}_N$ so that the upper boundary $\rho_N$ of $\mathcal{R}_N$ (the path on the dual hexagonal lattice bordering it above) is *a* front. We shall see a bit later that it is indeed *the* front (i.e., the only interface) with high probability.

We a priori expect it to be close to the line $\{y = 0\}$ corresponding to the sites where $p$ is critical. A hand-waving argument goes as follows: if we are at a distance approximately $N^{4/7}$ from the line, say, for instance, above it, the corresponding percolation parameter will be about $1/2 - N^{-3/7}$ and the associated characteristic length also of order $N^{4/7}$ [by using (11)]. Hence, the strip on which everything looks close to critical percolation has width of order $N^{4/7}$, and outside this strip, one is in the super- or subcritical regime. This intuitively explains the following result illustrated on Figure 5.

Theorem 6 (Localization).

- For all $\delta > 0$, there exists a $\delta' > 0$ such that for all sufficiently large $N$ and all $\ell_N \geq N^{4/7}$,

$$\mathbb{P}(\mathcal{R}_N \subseteq [\pm N^{4/7-\delta}]) \leq e^{-N^{\delta'}}. \tag{14}$$

- For all $\delta > 0$ and all $\gamma \geq 1$, there exists a $\delta' > 0$ such that for all sufficiently large $N$ and all $N^{4/7} \leq \ell_N \leq N^\gamma$,

$$\mathbb{P}(\mathcal{R}_N \not\subseteq [\pm N^{4/7+\delta}]) \leq e^{-N^{\delta'}}. \tag{15}$$

PROOF. Let us first derive the statement concerning $\mathbb{P}(\mathcal{R}_N \subseteq [\pm N^{4/7-\delta}])$. For that purpose, consider disjoint rhombi of the form

$$[i, i+2N^{4/7-\delta}] \times [-N^{4/7-\delta}, N^{4/7-\delta}] \qquad (i = 0, 2N^{4/7-\delta}+1, 4N^{4/7-\delta}+2, \ldots).$$

We can take at least $\ell_N/3N^{4/7-\delta}$ such rhombi and each of them possesses (independently of the others) a vertical vacant crossing with probability larger than

$$\mathbb{P}_{1/2+N^{-3/7-\delta}/2}(\mathcal{C}_H^*([0, 2N^{4/7-\delta}]^2)).$$

But [using (11), the critical exponent for $L$],

$$L(1/2 + N^{-3/7-\delta}/2) \approx (N^{-3/7-\delta}/2)^{-4/3} \approx N^{4/7+4\delta/3}$$

hence

$$L(1/2 + N^{-3/7-\delta}/2) \gg 2N^{4/7-\delta}$$

so that there exists a vertical vacant crossing with probability larger than $\varepsilon_0$. By independence, a "block" implying that $\mathcal{R}_N \not\subseteq [\pm N^{4/7-\delta}]$ will occur with probability larger than

$$1 - (1-\varepsilon_0)^{\ell_N/3N^{4/7-\delta}},$$

which proves the claim, as $\ell_N \geq N^{4/7}$ by assumption.

Let us now turn our attention to the quantity $\mathbb{P}(\mathcal{R}_N \not\subseteq [\pm N^{4/7+\delta}])$. Assume that $\mathcal{R}_N$ is not entirely contained in the strip $[\pm N^{4/7+\delta}]$.

If, at some point, it is, for example, above this strip, we are faced with the following alternatives:

- $\mathcal{R}_N$ is in the strip $[\pm N^{4/7+\delta}/2]$ at some point. In that case, consider the following rhombi, located between the lines $y = N^{4/7+\delta}/2$ and $y = N^{4/7+\delta}$: $[i, i+N^{4/7+\delta}/2] \times [N^{4/7+\delta}/2, N^{4/7+\delta}]$ $(i = 0, 1, 2, \ldots)$. It is easy to see (see Figure 6) that $\mathcal{R}_N$ will have to cross one of them vertically *or horizontally*.

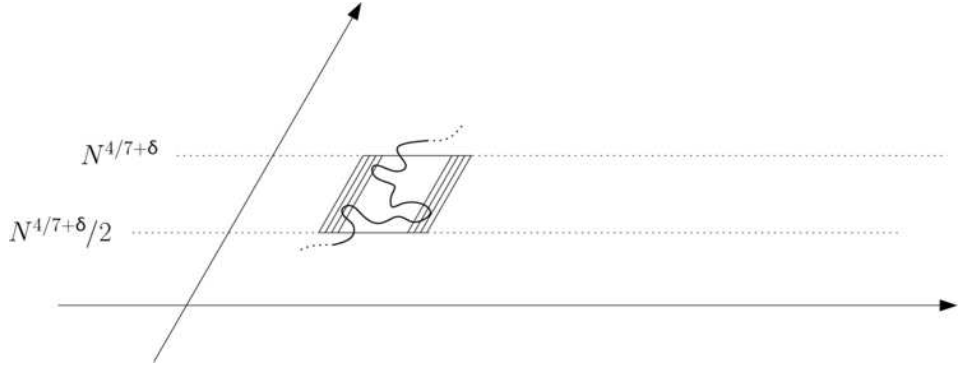

Fig. 6. $\mathcal{R}_N$ crosses one of the small rhombi vertically or horizontally.

There are at most $\ell_N$ such rhombi and they are in a zone for which $p \leq 1/2 - N^{-3/7+\delta}/4$. A crossing thus occurs with probability less than

$$\ell_N \mathbb{P}_{1/2-N^{-3/7+\delta}/4}(\mathcal{C}_H([0, N^{4/7+\delta}/2]^2) \cup \mathcal{C}_V([0, N^{4/7+\delta}/2]^2))$$

$$\leq 2\ell_N \mathbb{P}_{1/2-N^{-3/7+\delta}/4}(\mathcal{C}_H([0, N^{4/7+\delta}/2]^2))$$

$$\leq 2C\ell_N e^{-N^{4/7+\delta}/2L(1/2-N^{-3/7+\delta}/4)}$$

by Lemma 2 (of exponential decay).
But [by (11)],

$$L(1/2 - N^{-3/7+\delta}/4) \approx (N^{-3/7+\delta}/4)^{-4/3} \approx N^{4/7-4\delta/3}$$

and $\ell_N$ does not grow too fast, so the probability of the considered event tends to 0 subexponentially fast.
- $\mathcal{R}_N$ stays constantly above the strip $[\pm N^{4/7+\delta}/2]$, where the percolation parameter remains smaller than $1/2 - N^{-3/7+\delta}/4$. In that case, $\mathcal{R}_N$ will cross one of the rhombi $[0, N^{4/7+\delta}/2] \times [j + N^{4/7+\delta}/2, j + N^{4/7+\delta}]$ ($j = 0, 1, 2, \ldots$) forming a "column," vertically or horizontally.

There are at most $N$ such rhombi. Once again, Lemma 2 implies that a crossing occurs with probability less than

$$2CNe^{-N^{4/7+\delta}/2L(1/2-N^{-3/7+\delta}/4)}.$$

As before,

$$L(1/2 - N^{-3/7+\delta}/4) \approx N^{4/7-4\delta/3},$$

so the considered probability tends to 0 subexponentially fast.

If $\mathcal{R}_N$ is below the strip at some point, the argument is identical: consider instead the vacant crossing bordering $\mathcal{R}_N$. Hence, the final probability must simply be multiplied by 2. $\square$

3.3. *Uniqueness of the front.* In the previous subsection, we focused on the highest horizontal crossing $\mathcal{R}_N$. Clearly, the results remain valid if we instead consider the lowest horizontal vacant crossing $\mathcal{R}_N^*$. Recall that $\mathcal{R}_N$ is bordered above by a horizontal vacant crossing so that its upper boundary $\rho_N$ is a front. Similarly, $\mathcal{R}_N^*$ is bordered below by an occupied crossing and its lower boundary $\rho_N^*$ is also a front.

Note that $\rho_N^*$ is always below $\rho_N$. It is easy to see that confirming the uniqueness of the front amounts to checking that $\rho_N$ and $\rho_N^*$ coincide. It is also equivalent to verifying that $\mathcal{R}_N$ is connected to the bottom $\mathcal{B}_N$ by an occupied path (or that $\mathcal{R}_N^*$ is connected to the top $\mathcal{T}_N$ by a vacant path). We are now going to prove that this indeed occurs with a very large probability.

Note that if one starts an exploration process from the top-left corner of $\mathcal{S}_N$, one discovers the topmost occupied crossing without discovering the status of the sites below it. This RSW-type observation will be essential in our proof.

PROPOSITION 7 (Uniqueness). *Assume that $\ell_N \geq N^{4/7+\delta}$ for some $\delta > 0$. Then there is a $\delta' > 0$ (depending on $\delta$) such that for all sufficiently large $N$,*

$$\mathbb{P}(\rho_N = \rho_N^*) \geq 1 - e^{-N^{\delta'}}. \tag{16}$$

PROOF. As mentioned above, we will work with $\mathcal{R}_N$. Our goal is to show that the probability for $\mathcal{R}_N$ not to be connected to $\mathcal{B}_N$ is very small. For that purpose, we first divide the strip $\mathcal{S}_N$ into disjoint substrips $(\mathcal{S}_N^i)$, as follows. For $\varepsilon := \delta/4$, we choose $N^{3\varepsilon}/6$ disjoint substrips of length $3N^{4/7+\varepsilon}$ (not necessarily entirely covering $\mathcal{S}_N$) of the type

$$\mathcal{S}_N^i = [n_N^i, n_N^i + 3N^{4/7+\varepsilon}] \times [-N, N] \qquad (i = 0, \ldots, N^{3\varepsilon}/6 - 1).$$

Consider one of these strips $\mathcal{S}_N^i$. $\mathcal{R}_N$ crosses it horizontally and remains "below" its highest horizontal crossing, which we denote by $r_N^i$. Consequently, it will be sufficient to show that one of the $r_N^i$'s is connected to $\mathcal{B}_N$ in $\mathcal{S}_N^i$.

We now fix an $i$ and try to find a lower bound for the probability that $r_N^i$ is connected to the bottom of $\mathcal{S}_N^i$ by an occupied path that stays in that substrip. Let us suppose, for notational convenience, that $i = 0$ and $n_N^i = 0$. Note, first, that with probability at least $1/2$, there exists a vacant top-to-bottom crossing of the rhombus $[N^{4/7+\varepsilon}, 2N^{4/7+\varepsilon}] \times [0, N^{4/7+\varepsilon}]$ (percolation is subcritical in this region) so that a lowest point $z$ on $r_N^i$ in the middle part $[N^{4/7+\varepsilon}, 2N^{4/7+\varepsilon}] \times [-N, N]$ of the strip lies below the $x$-axis with probability at least $1/2 + o(1)$ [the localization result (Theorem 6) tells us that $r_N^i$ remains below the height $N^{4/7+\varepsilon}$ with high probability].

Now, we may have explored this highest crossing of $\mathcal{S}_N^0$ "from above," and have not yet discovered the status of the sites below it, so we can apply the FKG inequality for events involving only the state of these remaining sites. We are now going to show that in the case where $r_N^0$ passes below the $x$-axis, the conditional probability that it is connected to the bottom part of the substrip by an occupied crossing is bounded from below by a quantity of order $N^{-2\varepsilon}$.

A way to prove this is as follows. Let us first choose $z$ and define the annulus $S_{2N^{4/7-\varepsilon}}(z) \setminus \mathring{S}_{N^{4/7-\varepsilon}}(z)$ around $z$ depicted on Figure 7. Since it is contained in the region where $p \geq 1/2 - 2N^{-3/7-\varepsilon}$ and since the characteristic length corresponding to this value of the parameter is of order
$$L(1/2 - 2N^{-3/7-\varepsilon}) \approx N^{4/7+4\varepsilon/3},$$
there is a probability of at least $\delta_4^4$ (this is the constant coming from RSW theory) to observe an occupied circuit in this annulus.

We now want to connect this circuit to the bottom boundary of the substrip. Note that the part of the circuit that is below $r_N^0$, together with $r_N^0$, contains an occupied circuit around the segment $I = z + [-N^{4/7-\varepsilon}, N^{4/7-\varepsilon}] \times \{-N^{4/7-\varepsilon}\}$. We need the following simple lemma for critical percolation (see Figure 8).

LEMMA 8. *Consider the rhombus $[-N^{4/7+\varepsilon}, N^{4/7+\varepsilon}] \times [-2N^{4/7+\varepsilon}, 0]$ and the subinterval $I_N = [-N^{4/7-\varepsilon}, N^{4/7-\varepsilon}] \times \{0\}$ on its top edge. Then the event $\mathcal{C}_V^{I_N}([-N^{4/7+\varepsilon}, N^{4/7+\varepsilon}] \times [-2N^{4/7+\varepsilon}, 0])$ that there exists a vertical occupied crossing connecting $I_N$ to the bottom edge has a probability of at least*

(17) $$\mathbb{P}_{1/2}[\mathcal{C}_V^{I_N}([-N^{4/7+\varepsilon}, N^{4/7+\varepsilon}] \times [-2N^{4/7+\varepsilon}, 0])] \geq \frac{C}{N^{2\varepsilon}}$$

*for some universal constant $C$ (depending neither on $N$ nor on $\varepsilon$).*

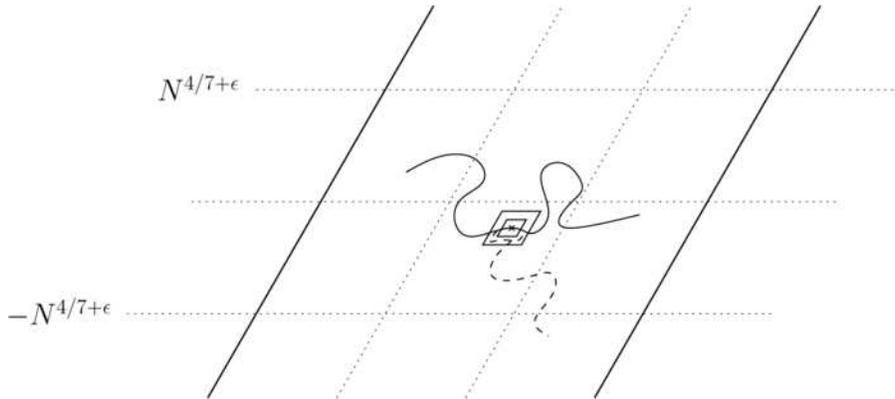

FIG. 7. *Construction of an arm going out of $[\pm N^{4/7+\varepsilon}]$.*

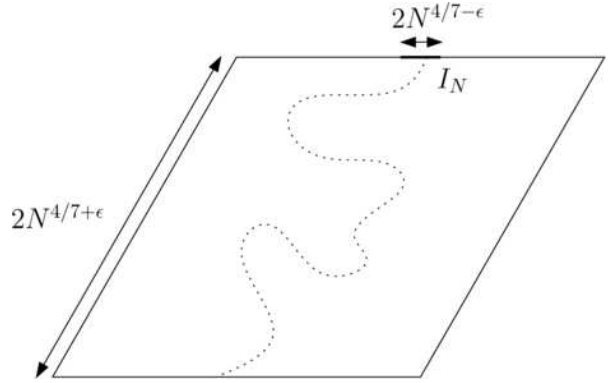

FIG. 8. *Existence of an arm from $I_N$ to the bottom side.*

PROOF. Consider the parallelogram $[0, N^{4/7+\varepsilon}] \times [-2N^{4/7+\varepsilon}, 0]$ and cover its top edge by less than $N^{2\varepsilon}$ intervals $I_N^j = [n_N'^j - N^{4/7-\varepsilon}, n_N'^j + N^{4/7-\varepsilon}] \times \{0\}$ of length $2N^{4/7-\varepsilon}$. We know from the RSW theorem that there exists a vertical occupied crossing with probability at least $\delta_2 > 0$, so

$$\delta_2 \leq \sum_j \mathbb{P}_{1/2}[\mathcal{C}_V^{I_N^j}([0, N^{4/7+\varepsilon}] \times [-2N^{4/7+\varepsilon}, 0])].$$

But, for each $j$,

$$\mathbb{P}_{1/2}[\mathcal{C}_V^{I_N^j}([0, N^{4/7+\varepsilon}] \times [-2N^{4/7+\varepsilon}, 0])]$$

$$\leq \mathbb{P}_{1/2}[\mathcal{C}_V^{I_N^j}([n_N'^j - N^{4/7+\varepsilon}, n_N'^j + N^{4/7+\varepsilon}] \times [-2N^{4/7+\varepsilon}, 0])]$$

$$= \mathbb{P}_{1/2}[\mathcal{C}_V^{I_N}([-N^{4/7+\varepsilon}, N^{4/7+\varepsilon}] \times [-2N^{4/7+\varepsilon}, 0])]$$

by translation invariance. Hence,

$$\delta_2 \leq N^{2\varepsilon} \times \mathbb{P}_{1/2}[\mathcal{C}_V^{I_N}([-N^{4/7+\varepsilon}, N^{4/7+\varepsilon}] \times [-2N^{4/7+\varepsilon}, 0])],$$

which completes the proof. Note that a repeated application of the RSW theorem [$\log(N^{2\varepsilon})$ times] would have given a lower bound of the type $N^{-\kappa\varepsilon}$ that would also have been sufficient for our purposes here. □

Putting the pieces together and making use of the FKG inequality, we get that for each $i \leq N^{3\varepsilon}/6 - 1$, the probability that $\mathcal{R}_N$ is connected to the bottom part of the strip in the substrip $\mathcal{S}_N^i$ is bounded from below by $C'/N^{2\varepsilon}$, independently for each $i$. The proposition then follows readily. Indeed, $\mathcal{R}_N$ is connected to $\mathcal{B}_N$ with probability at least

(18) $$1 - (1 - C'N^{-2\varepsilon})^{N^{3\varepsilon}/6} \geq 1 - e^{-N^{\varepsilon'}}$$

for some positive $\varepsilon'$. □

REMARK 9. On the other hand, note that if $\ell_N \leq N^{4/7-\delta}$ for some $\delta > 0$, then, by RSW, there exist several interfaces with probability bounded away from 0.

The previous results suggest that some restrictions should be made on the length $\ell_N$ of the strip. In the following, we will thus assume that there exists a $\delta > 0$ such that $\ell_N \geq N^{4/7+\delta}$: this hypothesis implies that the event corresponding to uniqueness has a probability tending to 1 subexponentially fast. For convenience, we also assume that $\ell_N = o(N^\gamma)$ for some $\gamma \geq 1$, to ensure that the front remains localized.

Hereafter, we will simply refer to *the* front and denote it by $\mathcal{F}_N$, which means that we will implicitly neglect the error term in the estimates that we derive. In particular, the front will be exactly the set of edges from which two arms can be drawn—one occupied to the bottom $\mathcal{B}_N$ and one vacant to the top $\mathcal{T}_N$.

**4. Length of the front.** We would now like to study the length $T_N$ of the front, that is, its number of edges. The preceding remark shows that we will need a two-arm probability estimate for that purpose. Unless otherwise stated, the expression "two arms" refers to two arms *of opposite color*.

We will often have to count edges rather than sites. Since we are primarily interested in rough estimates, we will then just use the fact the number of edges and the number of corresponding sites are comparable (i.e., up to a multiplicative factor of 6). To simplify notation, it will be convenient to associate to each edge $e$ one of its two neighboring sites, $x_e$, which we do arbitrarily and permanently.

4.1. *Two-arm estimates.* We are now going to derive the analog of (5) in the case of two arms and for nonconstant $p$. This lemma will enable us to estimate the probability of having two arms from an edge $e$ in the "critical strip." The goal is to show roughly that

$$\mathbb{P}_p[A^2(0,n)] \asymp \mathbb{P}_{1/2}[A^2(0,n)] \qquad (n \leq L(p)).$$

In fact, for our purposes, we will have to consider, instead of $\mathbb{P}_p$, product measures $P'$ with associated parameters $p'(v)$ which may depend on the site $v$, but remain between $p$ and $1-p$ (we will simply say that $P'$ is "between $\mathbb{P}_p$ and $\mathbb{P}_{1-p}$"). The present situation is a little more complicated than for the one-arm estimate because of the lack of monotonicity ($A^2$ events correspond to the combination of one path of each type, so they are neither increasing nor decreasing).

For a parallelogram $R$ and a site $v$ contained in its interior, Kesten considered in [14] the event $\Gamma(v, R)$ that there exist four arms, alternatively occupied and vacant, from the set $\partial v$ of vertices neighboring $v$ to the boundary

$\partial R$ of $R$:

$\Gamma(v, R) = \{$there exist four paths $r_1, \ldots, r_4$ from $\partial v$ to $\partial R$, ordered
clockwise, such that $r_1, r_3$ are occupied and $r_2, r_4$ are vacant$\}$.

Here, we will need this event in the course of proof, but we will be more interested in the analog for two arms:

$\Gamma_2(v, R) = \{$there exist an occupied path $r_1$ and a vacant path $r_2$ from $\partial v$ to $\partial R\}$.

Note that $A^2(0, n) = \Gamma_2(0, S_n)$. Let us now state and prove the result.

LEMMA 10. *Uniformly in $p$, $\hat{P}$ between $\mathbb{P}_p$ and $\mathbb{P}_{1-p}$, $n \leq L(p)$, we have*

(19) $$\hat{P}[\Gamma_2(0, S_n)] \asymp \mathbb{P}_{1/2}[\Gamma_2(0, S_n)].$$

PROOF. This proof is an adaptation of the proof of Theorem 1 in Kesten's paper [14]. We first recall some estimates contained in this paper and then adapt the original proof for one arm to the case of two arms. However, note that this result could be addressed as in [14]: Lemma 8 in this paper is the "four-arm" version of the present two-arm result.

Let us first introduce the events that we use throughout the proof. All of them are exact analogs of events defined in Kesten's paper [14]. We consider a parallelogram $R$ and a site $v$ such that $v \in \mathring{R}$.

- In the case of one arm, we do not lose much (a factor 4) by requiring the extremity to be on a specified edge of $R$, at least if $v$ is not too close to one of the edges. Here, we require the occupied arm to arrive at the bottom edge of $R$, and the vacant arm at the top edge. This event is denoted by $\Omega_2(v, R)$:

  $\Omega_2(v, R) = \{\Gamma_2(v, R)$ occurs, $r_1$ and $r_2$ having their extremities on the bottom and the top edges of $R$, respectively$\}$.

- The two paths $r_1$ and $r_2$ might be a bit difficult to extend, so we may want to add further "security strips." This leads to the definition of the event $\Delta_2$, analog in the case of two arms to the event $\Delta$. Although more restrictive, this event has a probability that remains comparable to $\Omega_2(0, R)$.

  We thus consider the two horizontal strips

  $$\mathcal{A}(1, k) := [-2^{k-1}, 2^{k-1}] \times [-2^k, -2^{k-1}],$$
  $$\mathcal{A}(2, k) := [-2^{k-1}, 2^{k-1}] \times [2^{k-1}, 2^k]$$

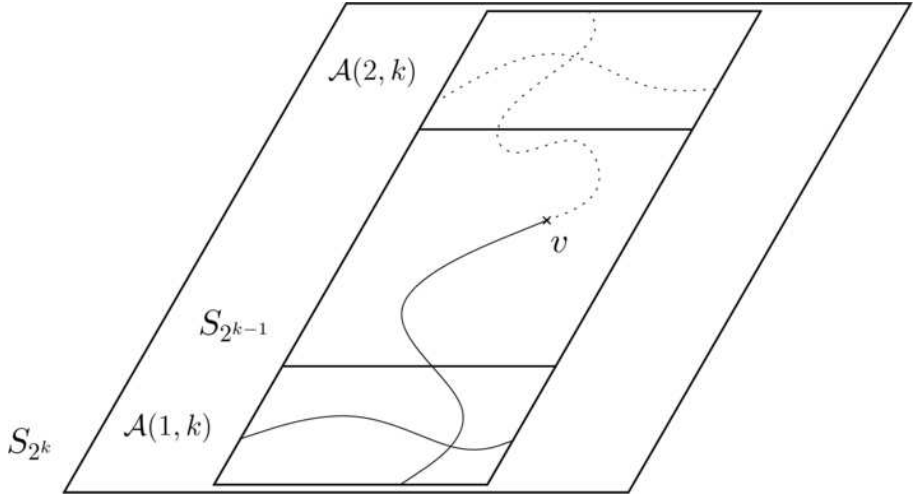

FIG. 9.  *Definition of the event* $\Delta_2(v, S_{2^k})$.

and define (see Figure 9), for a site $v$ in $S_{2^{k-1}}$,

$$\Delta_2(v, S_{2^k}) = \{\Gamma_2(v, S_{2^k}) \text{ occurs, with } r_i \cap (S_{2^k} \setminus \overset{\circ}{S}_{2^{k-1}}) \subseteq \mathcal{A}(i, k),$$
$$\text{and there exist two horizontal crossings,}$$
$$\text{one occupied of } \mathcal{A}(1, k) \text{ and one vacant of } \mathcal{A}(2, k)\}.$$

- We similarly define, for a parallelogram $R'$ contained in the interior of $S_{2^k}$,

$\tilde{\Gamma}_2(S_{2^k}, R') = \{$there exist an occupied path $r_1$ and a vacant path $r_2$ from $\partial R'$ to the bottom and top edges, respectively, of $S_{2^k}$ which are (with the exception of their extremities on $\partial R'$) contained in $S_{2^k} \setminus R'\}$

and for a site $v$ in $S_{2^{k-1}}$, $j \leq k-2$, the following strips (centered on $v$):

$$\tilde{\mathcal{A}}(1, j) := [v(1) - 2^j, v(1) + 2^j] \times [v(2) - 2^{j+1}, v(2) - 2^j],$$
$$\tilde{\mathcal{A}}(2, j) := [v(1) - 2^j, v(1) + 2^j] \times [v(2) + 2^j, v(2) + 2^{j+1}]$$

then

$$\tilde{\Delta}_2(S_{2^k}, S_{2^j}(v)) = \{\tilde{\Gamma}_2(S_{2^k}, S_{2^j}(v)) \text{ occurs, with associated paths } r_i$$
$$\text{satisfying } r_i \cap (S_{2^{j+1}}(v) \setminus \overset{\circ}{S}_{2^j}(v)) \subseteq \tilde{\mathcal{A}}(i, j),$$
$$\text{and there exist two horizontal crossings,}$$
$$\text{one occupied of } \tilde{\mathcal{A}}(1, j) \text{ and one vacant of } \tilde{\mathcal{A}}(2, j)\}.$$

Let us now state the estimates that we will use. We will not prove them since these are exact analogs for two arms of results stated for four arms in Kesten's paper [14], except that we generalize the condition "between $\mathbb{P}_{1/2}$ and $\mathbb{P}_p$" to "between $\mathbb{P}_p$ and $\mathbb{P}_{1-p}$." This generalization is valid since the only tool used in the proofs is the Russo–Seymour–Welsh theorem [we also implicitly use the symmetry about $1/2$, which implies that $L(p) = L(1-p)$]. In the following, $\hat{P}$ denotes any product measure "between $\mathbb{P}_p$ and $\mathbb{P}_{1-p}$."

(1) **Extension of $\Delta_2$.** There exists a constant $C_1 < \infty$ such that

(20) $$\hat{P}(\Delta_2(0, S_{2^k})) \leq C_1 \hat{P}(\Delta_2(0, S_{2^{k+1}}))$$

for all $p$, $\hat{P}$ and $2^k \leq L(p)$.

(2) **Comparison of $\Gamma_2$ and $\Delta_2$** (analog of Lemma 4 in [14]). There exists a constant $C_2 < \infty$ such that

(21) $$\hat{P}(\Gamma_2(0, S_{2^k})) \leq C_2 \hat{P}(\Delta_2(0, S_{2^k}))$$

for all $p$, $\hat{P}$ and $2^k \leq L(p)$.

(3) **Comparison of $\tilde{\Gamma}_2$ and $\tilde{\Delta}_2$** (analog of Lemma 5 in [14]). There exists a constant $C_3 < \infty$ such that

(22) $$\hat{P}(\tilde{\Gamma}_2(S_{2^k}, S_{2^j})) \leq C_3 \hat{P}(\tilde{\Delta}_2(S_{2^k}, S_{2^j}))$$

for all $p$, $\hat{P}$, $2^k \leq L(p)$ and $j \leq k-2$.

These prerequisites being recalled, we are now able to turn to the proof of Lemma 10. For that purpose, we will have to adapt Kesten's original proof. Recall that we consider a parameter $p$, a measure $\hat{P}$ between $\mathbb{P}_p$ and $\mathbb{P}_{1-p}$ and an integer $n \leq L(p)$. The parameters of $\hat{P}$ will be denoted by $p(v)$. We may assume that $p < 1/2$. Also, note that it suffices to prove the theorem with $\mathbb{P}_p$ instead of $\mathbb{P}_{1/2}$.

**1st step:** We first observe that we can replace $\Gamma_2(0, S_n)$ by $\Omega_2(0, S_{2^k})$, with $k$ such that $2^k \leq n < 2^{k+1}$. Indeed, the estimates (20) and (21) imply that

$$\hat{P}(\Gamma_2(0, S_n)) \geq \hat{P}(\Gamma_2(0, S_{2^{k+1}}))$$
$$\geq \hat{P}(\Delta_2(0, S_{2^{k+1}}))$$
$$\geq C_1^{-1} \hat{P}(\Delta_2(0, S_{2^k}))$$
$$\geq C_2^{-1} C_1^{-1} \hat{P}(\Gamma_2(0, S_{2^k}))$$
$$\geq C_2^{-1} C_1^{-1} \hat{P}(\Omega_2(0, S_{2^k}))$$

and

$$\hat{P}(\Omega_2(0, S_{2^k})) \geq \hat{P}(\Delta_2(0, S_{2^k}))$$
$$\geq C_2^{-1} \hat{P}(\Gamma_2(0, S_{2^k}))$$
$$\geq C_2^{-1} \hat{P}(\Gamma_2(0, S_n)).$$

**2nd step:** We would now like to use Russo's formula, but for technical reasons which will become clear in the next step, this will only work nicely for the points that are not too close to the boundary. We must therefore first see how the change of $p(v)$ for the points that are close to the boundary affects the probabilities that we investigate. More precisely, we start with $\mathbb{P}_p$ and change $p \rightsquigarrow p(v)$ only in $S_{2^k} \setminus S_{2^{k-3}}$. The resulting measure $\tilde{P}$ is between $\mathbb{P}_p$ and $\mathbb{P}_{1-p}$. Therefore,

$$\mathbb{P}_p(\Omega_2(0, S_{2^k})) \asymp \mathbb{P}_p(\Delta_2(0, S_{2^k})) \asymp \mathbb{P}_p(\Delta_2(0, S_{2^{k-3}}))$$

by (20) and (21), and, similarly,

$$\tilde{P}(\Omega_2(0, S_{2^k})) \asymp \tilde{P}(\Delta_2(0, S_{2^{k-3}})),$$

which allows us to conclude since $\mathbb{P}_p(\Delta_2(0, S_{2^{k-3}})) = \tilde{P}(\Delta_2(0, S_{2^{k-3}}))$ by definition.

**3rd step:** In order to handle the sites $v$ located in $S_{2^{k-3}}$, we will apply a generalization of Russo's formula to the family of measures $(\hat{P}_t)_{t \in [0,1]}$ with parameters

$$p(v, t) = tp(v) + (1-t)p,$$

which corresponds to a linear interpolation between $p \rightsquigarrow p(v)$ in $S_{2^{k-3}}$. The event $\Omega_2(0, S_{2^k})$ can be written as the intersection of $A_2(0, S_{2^k}) = \{$there exists an occupied path from $\partial 0$ to the bottom side of $S_{2^k}\}$ and $B_2(0, S_{2^k}) = \{$there exists a vacant path from $\partial 0$ to the top side of $S_{2^k}\}$. These two events are respectively increasing and decreasing and, in that case, $\frac{d}{dt}\hat{P}_t(A_2 \cap B_2)$ can be expressed as (see Lemma 1 in [14])

$$\sum_{v \in S_{2^{k-3}}} \frac{dp(v,t)}{dt}[\hat{P}_t(v \text{ is pivotal for } A_2, \text{ but not for } B_2, \text{ and } B_2 \text{ occurs})$$
$$- \hat{P}_t(v \text{ is pivotal for } B_2, \text{ but not for } A_2, \text{ and } A_2 \text{ occurs})].$$

A vertex $v \in S_{2^{k-3}}$, $v \neq 0$, is pivotal for $A_2$ iff there exist two paths, both containing $v$, such that:

1. their sites, except for $v$, are respectively all occupied and all vacant;
2. the first path connects $\partial 0$ to the bottom side of $S_{2^k}$;
3. the second path separates 0 from the bottom side of $S_{2^k}$.

In the case of one arm, Kesten used a vacant loop around the origin; here, the separating path can have its extremities on the boundary. This will not change the computations since we still have four arms locally that we will sum in the same way.

Actually, we must also assume that $k \geq 7$ and put aside the vertices which are too close to the origin, for instance, those in $S_{16}$. For these sites,

$$\hat{P}_t(v \text{ is pivotal for } A_2, \text{ but not for } B_2, \text{ and } B_2 \text{ occurs})$$
$$\leq \hat{P}_t(\tilde{\Gamma}_2(S_{2^k}, S_{16}))$$
$$\leq C_4 \hat{P}_t(\Omega_2(0, S_{2^k}))$$

for some universal constant $C_4$.

We now associate to each site $v \in S_{2^{k-3}} \setminus S_{16}$ a parallelogram $R(v)$ such that $0 \notin R(v)$ so that we make appear four arms locally. We choose them to present the following property: if $j$ is such that $2^{j+1} < d(v,0) \leq 2^{j+2}$, then $R(v)$ is included in $\overset{\circ}{S}_{2^{j+3}} \setminus S_{2^j}$. This will imply that

$$\hat{P}_t(v \text{ is pivotal for } A_2, \text{ but not for } B_2, \text{ and } B_2 \text{ occurs})$$
$$\leq \hat{P}_t[\Gamma_2(0, S_{2^j}) \cap \Gamma(v, R(v)) \cap \tilde{\Gamma}_2(S_{2^k}, S_{2^{j+3}})]$$
$$= \hat{P}_t[\Gamma_2(0, S_{2^j})]\hat{P}_t[\Gamma(v, R(v))]\hat{P}_t[\tilde{\Gamma}_2(S_{2^k}, S_{2^{j+3}})],$$

by independence (these events depending on sites in disjoint sets).

We can take the parallelograms $R(v)$ as in Kesten's paper (see [14] page 117). As observed in this paper, the precise choice is not really important. On one hand, we must ensure that the four arms are not too small. For that purpose, the distances between $v$ and each of the sides of $R(v)$ must be comparable to the distance between $0$ and $v$. On the other hand, if the announced property is satisfied, we will be in a position to join the paths outside $R(v)$ that are respectively between $0$ and $\partial S_{2^j}$ and between $\partial S_{2^{j+3}}$ and $\partial S_{2^k}$. For the sake of completeness, let us briefly recall how Kesten chooses these $R(v)$. Consider $v = (v_1, v_2) \notin S_{16}$. If $|v_1| \leq |v_2| \leq 2^{k-3}$ and $16 \leq 2^{j+1} < v_2 \leq 2^{j+2}$, take $l_1$ and $l_2$ such that

$$l_1 2^{j-2} < v_1 \leq (l_1+1)2^{j-2} \quad \text{and} \quad l_2 2^{j-2} < v_2 \leq (l_2+1)2^{j-2},$$

and define

$$R(v) = [(l_1-2)2^{j-2}, (l_1+2)2^{j-2}] \times [l_2 2^{j-2} - 2^j, l_2 2^{j-2} + 2^j].$$

If $v_2 < 0$, take the image of $R((v_1, -v_2))$ under the symmetry with respect to the $x$-axis. Finally, if $|v_2| < |v_1|$, simply exchange the roles of the first and second coordinates. We can easily check that these parallelograms possess, by construction, the desired property.

Now, by combining the two terms $\hat{P}_t[\Gamma_2(0, S_{2^j})]$ and $\hat{P}_t[\tilde{\Gamma}_2(S_{2^k}, S_{2^{j+3}})]$, we can produce $\hat{P}_t[\Omega_2(0, S_{2^k})]$: indeed, we easily get, from the Russo–Seymour–Welsh theorem (see Lemma 6 in [14]),

$$\delta_{16}^2 \hat{P}_t[\Delta_2(0, S_{2^j})]\hat{P}_t[\tilde{\Delta}_2(S_{2^k}, S_{2^{j+3}})] \leq \hat{P}_t[\Omega_2(0, S_{2^k})]$$

as, by assumption, $2^k \leq L(p)$ and $j \leq k - 5$ (here, a slight generalization of the FKG inequality is needed (see, e.g., Lemma 3 in [14]), invoking zones where "everything is monotonic").

We then obtain from

$$\hat{P}_t[\Delta_2(0, S_{2^j})] \geq C_2^{-1} \hat{P}_t[\Gamma_2(0, S_{2^j})]$$

[which follows from (21)] and

$$\hat{P}_t[\tilde{\Delta}_2(S_{2^k}, S_{2^{j+3}})] \geq C_3^{-1} \hat{P}_t[\tilde{\Gamma}_2(S_{2^k}, S_{2^{j+3}})]$$

[using (22)] that

(23) $$\hat{P}_t[\Gamma_2(0, S_{2^j})]\hat{P}_t[\tilde{\Gamma}_2(S_{2^k}, S_{2^{j+3}})] \leq C_5 \hat{P}_t[\Omega_2(0, S_{2^k})].$$

Hence,

$$\hat{P}_t(v \text{ is pivotal for } A_2, \text{ but not for } B_2, \text{ and } B_2 \text{ occurs})$$
$$\leq C_5 \hat{P}_t[\Omega_2(0, S_{2^k})]\hat{P}_t[\Gamma(v, R(v))].$$

For reasons of symmetry, the same is true if we invert $A_2$ and $B_2$. Consequently, dividing by $\hat{P}_t[\Omega_2(0, S_{2^k})]$ will produce its logarithmic derivative in the left-hand side of Russo's formula:

$$\left|\frac{d}{dt} \log[\hat{P}_t(\Omega_2(0, S_{2^k}))]\right|$$
$$\leq 2C_5 \sum_{v \in S_{2^{k-3}} \setminus S_{16}} \left|\frac{dp(v,t)}{dt}\right| \hat{P}_t[\Gamma(v, R(v))] + 2C_4 \sum_{v \in S_{16}} \left|\frac{dp(v,t)}{dt}\right|.$$

Finally, the first term that we obtain (the sum corresponding to the existence of four arms locally) is exactly the same as in Kesten's paper [14] (end of the proof of Theorem 1, page 140) and its integral between 0 and 1 is bounded by some universal constant $C_6$:

$$\int_0^1 \left(\sum_{v \in S_{2^{k-3}} \setminus S_{16}} \left|\frac{dp(v,t)}{dt}\right| \hat{P}_t[\Gamma(v, R(v))]\right) dt \leq C_6.$$

The desired conclusion then follows. $\square$

4.2. *Expected length of the front.* We are now able to study properties of the front in the "critical strip." Roughly speaking, only the sites in $[\pm N^{4/7}]$ must be taken into account and each of these sites has a probability of approximately $(N^{4/7})^{-1/4} = N^{-1/7}$ (the two-arm exponent being equal to $1/4$) to be on the front. Starting from this idea, we will prove the following estimate on the expectation of $T_N$.

CRITICAL EXPONENTS OF PLANAR GRADIENT PERCOLATION 23

PROPOSITION 11. *Recall that, by assumption, $\ell_N \geq N^{4/7+\delta}$ and $\ell_N = o(N^\gamma)$. For all $\varepsilon > 0$, we have, for $N$ sufficiently large,*

$$(24) \qquad N^{3/7-\varepsilon}\ell_N \leq \mathbb{E}[T_N] \leq N^{3/7+\varepsilon}\ell_N.$$

PROOF. Throughout the proof, we will use the fact that

$$(25) \qquad \mathbb{E}[T_N] = \sum_{e \in \mathcal{S}_N} \mathbb{P}(e \in \mathcal{F}_N).$$

We first consider the upper bound. Taking $\varepsilon' = \varepsilon/4$, we have

$$\mathbb{E}[T_N] \leq 6|\mathcal{S}_N| \times \mathbb{P}(\mathcal{F}_N \not\subseteq [\pm N^{4/7+\varepsilon'}]) + \sum_{e \in [\pm N^{4/7+\varepsilon'}]} \mathbb{P}(e \in \mathcal{F}_N)$$

and it follows from localization and the fact that $|\mathcal{S}_N| = (2N+1) \times (\ell_N+1) = o(N^{\gamma+1})$ that the first term tends to 0 subexponentially fast. We can thus restrict the summation to the vertices in the strip $[\pm N^{4/7+\varepsilon'}]$.

But, for $e \in [\pm N^{4/7+\varepsilon'}]$,

$$\mathbb{P}(e \in \mathcal{F}_N) \leq \mathbb{P}[2 \text{ arms } x_e \rightsquigarrow \partial S_{N^{4/7-2\varepsilon'}}(x_e)]$$

and as $S_{N^{4/7-2\varepsilon'}}(x_e) \subseteq [\pm 2N^{4/7+\varepsilon'}]$, the percolation parameter in this box remains in the range $[1/2 \pm 2N^{-3/7+\varepsilon'}]$. The associated characteristic length being

$$L(1/2 \pm 2N^{-3/7+\varepsilon'}) \approx N^{4/7-4\varepsilon'/3} \gg N^{4/7-2\varepsilon'},$$

we get (using Lemma 10)

$$\mathbb{P}[2 \text{ arms } x_e \rightsquigarrow \partial S_{N^{4/7-2\varepsilon'}}(x_e)]$$
$$\asymp \mathbb{P}_{1/2}[2 \text{ arms } x_e \rightsquigarrow \partial S_{N^{4/7-2\varepsilon'}}(x_e)]$$
$$\approx (N^{4/7-2\varepsilon'})^{-1/4} \qquad \text{(by using the 2-arm exponent)}$$
$$\ll N^{-1/7+\varepsilon'}.$$

Now we must simply sum this inequality to get the desired result: for $N$ large enough,

$$\sum_{e \in [\pm N^{4/7+\varepsilon'}]} \mathbb{P}(e \in \mathcal{F}_N) \leq 6(2N^{4/7+\varepsilon'} + 1)(\ell_N + 1) \times N^{-1/7+\varepsilon'}$$
$$\leq N^{3/7+\varepsilon}\ell_N.$$

Let us now turn to the lower bound. We restrict ourselves to the edges $e$ in the strip $[2N^{4/7+\varepsilon'}, \ell_N - 2N^{4/7+\varepsilon'}] \times [\pm N^{4/7-\varepsilon'}]$, with $\varepsilon' = \varepsilon/6$. For such edges, we would like to estimate the probability of having two arms, one occupied to $\mathcal{B}_N$ and one vacant to $\mathcal{T}_N$, so we will use the event $\Delta_2$ rather than $\Gamma_2$. Indeed, take $j$ such that $N^{4/7-\varepsilon'} < 2^j \leq 2N^{4/7-\varepsilon'}$: the probability of having two arms, one occupied to the bottom of $\partial S_{2^j}(x_e)$ and one vacant to the top, is at least

$$\mathbb{P}(\Delta_2(x_e, S_{2^j}(x_e))).$$

These two paths can then be extended so that they go out of the strip $[\pm N^{4/7+\varepsilon'}]$: Lemma 8 (contained in the proof of uniqueness) implies that this can be done with probability at least

(26) $$(CN^{-2\varepsilon'})^2 = C'N^{-4\varepsilon'}.$$

On the other hand, since we stay in the strip $[\pm 3N^{4/7-\varepsilon'}]$ of associated characteristic length $L(1/2 \pm 3N^{-3/7-\varepsilon'}) \approx N^{4/7+4\varepsilon'/3}$, we get that

$$\mathbb{P}(\Delta_2(x_e, S_{2^j}(x_e))) \asymp \mathbb{P}(\Gamma_2(x_e, S_{2^j}(x_e)))$$
$$\asymp \mathbb{P}_{1/2}(\Gamma_2(x_e, S_{2^j}(x_e)))$$
$$\geq \mathbb{P}_{1/2}(\Gamma_2(0, S_{2N^{4/7-\varepsilon'}}))$$

and the 2-arm exponent implies that

(27) $$\mathbb{P}_{1/2}(\Gamma_2(0, S_{2N^{4/7-\varepsilon'}})) \approx (2N^{4/7-\varepsilon'})^{-1/4} \gg N^{-1/7}.$$

We can thus construct two arms going out of the strip $[\pm N^{4/7+\varepsilon'}]$ with probability at least $N^{-1/7-4\varepsilon'}$ (for $N$ large enough).

In that case, $e$ has a high probability of being connected to the top and to the bottom of $\mathcal{S}_N$. Indeed, the front would go out of the strip $[\pm N^{4/7+\varepsilon'}]$ otherwise, which occurs with a probability less that $\varepsilon_N$, for some $\varepsilon_N$ (independent of $e$) tending to 0 subexponentially fast. The conclusion follows by summing the lower bound over all edges $e$ in the strip $[2N^{4/7+\varepsilon'}, \ell_N - 2N^{4/7+\varepsilon'}] \times [\pm N^{4/7-\varepsilon'}]$: for $N$ large enough,

$$\sum_{e \in \mathcal{S}_N} \mathbb{P}(e \in \mathcal{F}_N) \geq (2N^{4/7-\varepsilon'})(\ell_N - 4N^{4/7+\varepsilon'}) \times (N^{-1/7-4\varepsilon'} - \varepsilon_N)$$
$$\geq N^{3/7-\varepsilon}\ell_N. \qquad \square$$

4.3. *Convergence in $L^2$.* As seen in Section 4.2, everything happens in the strip $[\pm N^{4/7}]$ and it is possible to determine whether or not an edge $e$ is on the front on a neighborhood of size $N^{4/7}$ (with probability very close to 1). Distant points will thus be almost completely decorrelated and a phenomenon of horizontal "averaging" therefore occurs so that we get the following bound on $\mathrm{Var}[T_N]$.

PROPOSITION 12. *We still assume that $\ell_N \geq N^{4/7+\delta}$ and $\ell_N = o(N^\gamma)$. For each $\varepsilon > 0$, we have, for $N$ sufficiently large,*

(28) $$\mathrm{Var}[T_N] \leq N^{10/7+\varepsilon}\ell_N.$$

By combining this result with the estimates for $\mathbb{E}[T_N]$ from the previous subsection, we immediately get that

$$\text{Var}[T_N] = o(\mathbb{E}[T_N]^2) \tag{29}$$

and, consequently, the following.

THEOREM 13. *If $\ell_N \geq N^{4/7+\delta}$ and $\ell_N = o(N^\gamma)$, then*

$$\frac{T_N}{\mathbb{E}[T_N]} \longrightarrow 1 \quad \text{in } L^2 \text{ as } N \to \infty. \tag{30}$$

REMARK 14. On each vertical line, among the (approximately) $N^{4/7}$ edges which lie in the critical strip, roughly $N^{3/7}$ of them will be on the front, that is, a proportion of $1/N^{1/7}$. If we take $\ell_N$ to be equal to $N$, we get $\mathbb{E}[T_N] \sim N^{10/7}$ and this exponent can indeed be observed numerically.

PROOF OF PROPOSITION 12. For an edge $e$, the event "$e \in \mathcal{F}_N$" will be denoted by $F_e$. Using the expression

$$T_N = \sum_e \mathbb{I}_{F_e}, \tag{31}$$

we get

$$\text{Var}[T_N] = \sum_{e,f}[\mathbb{P}(F_e \cap F_f) - \mathbb{P}(F_e)\mathbb{P}(F_f)]. \tag{32}$$

First, take $\varepsilon' = \varepsilon/8$ and note that we can restrict the summation to the edges $e, f \in [\pm N^{4/7+\varepsilon'}]$: similarly to the argument following (25), the remaining term tends to 0 subexponentially fast.

The idea is to replace $F_e$ by an event $\tilde{F}_e$ which depends only on sites in a box around $e$ of size $N^{4/7+\varepsilon'}$ and such that $\mathbb{P}(F_e \Delta \tilde{F}_e) \to 0$ subexponentially fast (uniformly in $e$). To construct such an event, we simply invoke the result of uniqueness in the strip of length $2N^{4/7+\varepsilon'}$ centered horizontally on $x_e$: it implies that we can take

$$\tilde{F}_e := \{2 \text{ arms } e \rightsquigarrow \text{ top and bottom sides of } \partial(S_{N^{4/7+\varepsilon'}}(x_e) \cap \mathcal{S}_N)\}.$$

We thus have

$$\text{Var}[T_N] = \sum_e \sum_f [\mathbb{P}(\tilde{F}_e \cap \tilde{F}_f) - \mathbb{P}(\tilde{F}_e)\mathbb{P}(\tilde{F}_f)] + \varepsilon_N$$

with an error term $\varepsilon_N$ tending to 0 subexponentially fast.

We now fix an edge $e$ and estimate the corresponding sum. Only the edges $f$ in $S_{2N^{4/7+\varepsilon'}}(x_e)$ must be taken into account; otherwise, the two associated boxes do not intersect and

$$\mathbb{P}(\tilde{F}_e \cap \tilde{F}_f) - \mathbb{P}(\tilde{F}_e)\mathbb{P}(\tilde{F}_f) = 0.$$

We now have (and it is enough) to estimate
$$\sum_{f \in S_{2N^{4/7+\varepsilon'}}(x_e) \cap \mathcal{S}_N} \mathbb{P}(\tilde{F}_e \cap \tilde{F}_f).$$

Let $d = d(x_e, x_f)$. If $d > N^{4/7-2\varepsilon'}/4$, then we have
$$\mathbb{P}(\tilde{F}_e \cap \tilde{F}_f) \leq \mathbb{P}(\Gamma_2(x_e, S_{N^{4/7-2\varepsilon'}/8}(x_e)) \cap \Gamma_2(x_f, S_{N^{4/7-2\varepsilon'}/8}(x_f)))$$
$$= \mathbb{P}(\Gamma_2(x_e, S_{N^{4/7-2\varepsilon'}/8}(x_e)))\mathbb{P}(\Gamma_2(x_f, S_{N^{4/7-2\varepsilon'}/8}(x_f)))$$
$$\leq N^{-2/7+2\varepsilon'},$$

for $N$ large enough, by using the 2-arm exponent in the boxes $S_{N^{4/7-2\varepsilon'}/8}(x_e)$ and $S_{N^{4/7-2\varepsilon'}/8}(x_f)$. The sum over the edges $f$ in $S_{2N^{4/7+\varepsilon'}}(x_e) \setminus S_{N^{4/7-2\varepsilon'}/4}(x_e)$ is thus at most
$$6(4N^{4/7+\varepsilon'} + 1)^2 N^{-2/7+2\varepsilon'} \leq N^{6/7+5\varepsilon'}.$$

Now, assume that $4 \leq d \leq N^{4/7-2\varepsilon'}/4$. By introducing two boxes of size $d/2$ around $x_e$ and $x_f$, and another box of size $2d$ around $x_e$, we see that

$$\mathbb{P}(\tilde{F}_e \cap \tilde{F}_f)$$
$$\leq \mathbb{P}(\Gamma_2(x_e, S_{d/2}(x_e)))\mathbb{P}(\Gamma_2(x_f, S_{d/2}(x_f)))\mathbb{P}(\tilde{\Gamma}_2(S_{2N^{4/7-2\varepsilon'}}(x_e), S_{2d}(x_e))).$$

(Here, to be completely rigorous, we should prove that the condition on the extremities of the arms in $\tilde{\Gamma}_2$ can be relaxed: the methods of [14] apply—see also [18].) We can then combine the first and the third terms, as in (23):
$$\mathbb{P}(\Gamma_2(x_e, S_{d/2}(x_e)))\mathbb{P}(\tilde{\Gamma}_2(S_{2N^{4/7-2\varepsilon'}}(x_e), S_{2d}(x_e)))$$
$$\leq C_1 \mathbb{P}(\Gamma_2(x_e, S_{2N^{4/7-2\varepsilon'}}(x_e)))$$
$$\leq N^{-1/7+\varepsilon'}.$$

Hence, using the fact that there are at most $C_2 j$ edges at a distance $j$ from $e$ and the 2-arm exponent,
$$\sum_{f \in (S_{N^{4/7-2\varepsilon'}/4}(x_e) \setminus S_4(x_e)) \cap \mathcal{S}_N} \mathbb{P}(\tilde{F}_e \cap \tilde{F}_f)$$
$$\leq \sum_{j=4}^{N^{4/7-2\varepsilon'}/4} (C_2 j) \times (N^{-1/7+\varepsilon'} \mathbb{P}(\Gamma_2(x_f, S_{j/2}(x_f))))$$
$$\leq C_3 N^{-1/7+\varepsilon'} \sum_{j=4}^{N^{4/7-2\varepsilon'}/4} j \times j^{-1/4+\varepsilon'}$$
$$\leq C_4 N^{-1/7+\varepsilon'} (N^{4/7-2\varepsilon'}/4)^{7/4+\varepsilon'}$$
$$\leq N^{6/7+2\varepsilon'}.$$

Finally, the edges at a distance $d < 4$ make a contribution which is at most $C_5 \mathbb{P}(\tilde{F}_e) \leq N^{-1/7+\varepsilon'} \ll N^{6/7+2\varepsilon'}$.

Summing the different contributions, we get that

$$\sum_{f \in S_{2N^{4/7+\varepsilon'}}(x_e) \cap \mathcal{S}_N} \mathbb{P}(\tilde{F}_e \cap \tilde{F}_f) \leq N^{6/7+6\varepsilon'}.$$

Hence,

$$\begin{aligned} \operatorname{Var}[T_N] &\leq [6(2N^{4/7+\varepsilon'} + 1)(\ell_N + 1)] \times [N^{6/7+6\varepsilon'}] + \varepsilon_N \\ &\leq N^{10/7+\varepsilon} \ell_N. \end{aligned}$$  □

REMARK 15. Heuristically, this bound can be explained as follows: if we divide the strip $\mathcal{S}_N$ into approximately $(\ell_N/N^{4/7})$ disjoint (and thus more or less independent) boxes, the contribution of each of them to the variance is of order $((N^{4/7})^2)^2 N^{-2/7} = N^{14/7}$ so that the total variance is of order

$$(\ell_N/N^{4/7}) \times N^{14/7} = N^{10/7} \ell_N.$$

**5. Outer boundary.** To further describe $\mathcal{F}_N$ (still assuming it is unique), we can also consider its outer boundary (or accessible perimeter). Actually, two curves arise: the "upper" outer boundary $\mathcal{U}_N^+$ and the "lower" outer boundary $\mathcal{U}_N^-$. The upper outer boundary can be defined, for instance, as the lowest *self-avoiding* vacant crossing of $\mathcal{S}_N$. It can be viewed as the upper boundary of $\mathcal{F}_N$ "without the fjords" (the sites that are connected to the top of $\mathcal{S}_N$ only by a one-site-wide passage). For a site $v$, being on $\mathcal{U}_N^+$ is equivalent to have three arms: two disjoint vacant arms to the left and right sides of $\mathcal{S}_N$ and an occupied arm to the bottom.

The lengths of $\mathcal{U}_N^+$ and $\mathcal{U}_N^-$, denoted respectively by $U_N^+$ and $U_N^-$, can be computed in a similar way as the length $T_N$ of the front, by counting the 3-arm sites (two vacant and one occupied for $U_N^+$, two occupied and one vacant for $U_N^-$). There is, however, an extra technical difficulty involved in adapting Kesten's result: being pivotal still leads to four arms locally, but not necessarily long enough (due to consecutive arms of the same color) and the reasoning must be adapted (see [18] for more details). We get the following estimate on $\mathbb{E}[U_N^\pm]$.

PROPOSITION 16. *We still assume that $\ell_N \geq N^{4/7+\delta}$ and $\ell_N = o(N^\gamma)$. For all $\varepsilon > 0$, we have, for $N$ sufficiently large,*

(33) $$N^{4/21-\varepsilon} \ell_N \leq \mathbb{E}[U_N^\pm] \leq N^{4/21+\varepsilon} \ell_N.$$

The upper bound can be treated in the same way as in Proposition 11, but for the lower bound, we must be a little more careful: we must ensure that we can extend the two arms of the same type into two *disjoint* macroscopic arms. It can be done by using the fact that we can impose an occupied crossing in an $N^{4/7+\varepsilon} \times N^{4/7+\varepsilon}$ rhombus to arrive in a "corner" of size $N^{4/7-\varepsilon}$ with probability at least $N^{-\kappa\varepsilon}$ for some constant $\kappa$ (e.g., applying the RSW estimate $\log(N^{2\varepsilon})$ times; see Figure 10).

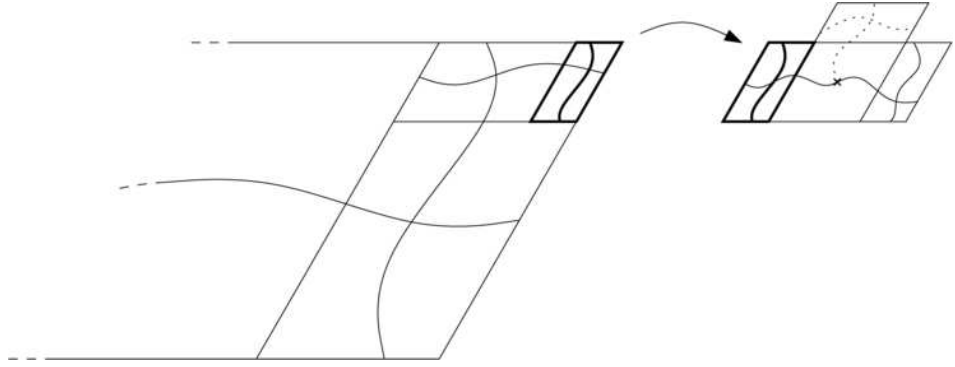

Fig. 10. *We can extend the three arms by successive applications of RSW.*

The decorrelation arguments that we used for the length of the front can also be applied to the outer boundary, implying convergence in $L^2$.

PROPOSITION 17. *If $\ell_N \geq N^{4/7+\delta}$ and $\ell_N = o(N^\gamma)$, then for each $\varepsilon > 0$, we have, for $N$ sufficiently large,*

$$\text{Var}[U_N^\pm] \leq N^{20/21+\varepsilon} \ell_N. \tag{34}$$

*Hence,*

$$\frac{U_N^\pm}{\mathbb{E}[U_N^\pm]} \longrightarrow 1 \qquad \text{in } L^2 \text{ as } N \to \infty. \tag{35}$$

Note that when $\ell_N = N$, this result implies that $U_N^\pm$ is of order $N^{25/21}$ with high probability.

**6. Open questions.** We have here used the logarithmic equivalence of Proposition 1 that was proven in [23]. It is, in fact, conjectured that this equivalence holds up to multiplicative constants. For instance, the stronger assumption $\mathbb{P}_{1/2}(A_n^4) \asymp n^{-5/4}$ would lead to $L(p) \asymp |p - 1/2|^{-4/3}$, which would make it possible to derive a sharper description of the front.

A natural issue would be to describe (if it exists) the scaling limit of the front, once properly renormalized. We could first argue as follows: in each of the substrips $[\pm N^{4/7-\varepsilon}]$, everything "looks like criticality" so that some aspects of the front will be the same as those of critical percolation interfaces, that is, $SLE_6$ curves.

However, we do not know at present how the front "bounces," that is to say, what happens when we are at a distance of approximately $N^{4/7}$ from the critical line. In fact, the front could be linked to the objects Camia, Fontes and Newman have recently introduced to describe percolation near criticality ([6, 7]) when we take the parameter to be $p = 1/2 + \lambda \delta^{3/4}$ with a possibly inhomogeneous $\lambda$.

**Acknowledgments.** The author would like to warmly thank his supervisor W. Werner for suggesting this problem to him and for his guidance. He also wishes to thank J. van den Berg and F. Camia for useful discussions and remarks.

DÉPARTEMENT DE MATHÉMATIQUES ET APPLICATIONS
ÉCOLE NORMALE SUPÉRIEURE
45 RUE D'ULM
75230 PARIS CEDEX 05
FRANCE
E-MAIL: pierre.nolin@ens.fr
AND
LABORATOIRE DE MATHÉMATIQUES
BÂTIMENT 425
UNIVERSITÉ PARIS-SUD
91405 ORSAY CEDEX
FRANCE